\documentclass{article}

\usepackage{amsmath,amsthm,amssymb,latexsym}
\usepackage{graphicx}
\usepackage{subfigure}

\newtheorem{thm}{Theorem}
\newtheorem{lem}[thm]{Lemma}

\theoremstyle{definition}
\newtheorem{dfn}[thm]{Definition}

\theoremstyle{remark}
\newtheorem*{rmk}{Remark}
\newtheorem*{rmks}{Remarks}

\newcommand{\crn}{\mathrm{cr}}

\title{Simultaneous Arithmetic Progressions on Algebraic Curves}
\author{Ryan Schwartz \and J\'ozsef Solymosi\footnote{The second author was supported by a Sloan Fellowship and NSERC and OTKA  grants.} \and Frank de Zeeuw \\ \\ Department of Mathematics \\ University of British Columbia \\ Vancouver, B.C., Canada V6T1Z2 \\ Email: \{ryano,solymosi,fdezeeuw\}@math.ubc.ca}

\begin{document}
\maketitle

\begin{abstract}
  A \emph{simultaneous arithmetic progression} (s.a.p.) of length $k$
  consists of $k$ points $(x_i, y_{\sigma(i)})$, where $x_i$ and $y_i$
  are arithmetic progressions and $\sigma$ is a permutation.
  Garcia-Selfa and Tornero asked whether there is a bound on the
  length of an s.a.p.~on an elliptic curve in Weierstrass form over
  $\mathbb{Q}$.  We show that $4319$ is such a bound for curves over
  $\mathbb{R}$.  This is done by considering translates of the curve
  in a grid as a graph.  A simple upper bound is found for the number
  of crossings and the ``crossing inequality'' gives a lower bound.
  Together these bound the length of an s.a.p.~on the curve.  We then
  use a similar method to extend the result to arbitrary real
  algebraic curves.  Instead of considering s.a.p.'s we consider
  $k^{2/3}$ points in a grid.  The number of crossings is bounded by
  B\'ezout's Theorem.  We then give another proof using a result of
  Jarn\'{\i}k bounding the number of grid points on a convex curve.
  This result applies as any real algebraic curve can be broken up
  into convex and concave parts, the number of which depend on the
  degree.  Lastly, these results are extended to complex algebraic
  curves.
\end{abstract}

\section{Introduction}

There are interesting problems in number theory related to arithmetic
progressions on elliptic curves. An example of such an open problem
is, what is the maximum number (if such a number exists) of rational
points on an elliptic curve such that their x-coordinates are in
arithmetic progression?  In \cite{Brem99b}, Bremner found elliptic
curves in Weierstrass form with arithmetic progressions of length 8 on
them, and Campbell found elliptic curves of the form $y^2 = f(x)$, with
$f$ a quartic, that contain arithmetic progressions of length $12$. In
\cite{Brem99a}, Bremner described how these arithmetic progressions
are related to $3\times 3$ magic squares with square
entries. Silverman, Bremner and Tzanakis noted in \cite{BrST00} that points in
arithmetic progression on elliptic curves are often independent with
respect to the group structure, which suggests a relation with the
much-researched rank of the curve.

In \cite{Garc06}, Garcia-Selfa and Tornero looked instead for
``simultaneous'' arithmetic progressions on elliptic curves, which are
defined as follows.
\begin{dfn}
  A \emph{simultaneous arithmetic progression} (s.a.p.) of length $k$ consists of 
  points $(x_i, y_{\sigma(i)})$, where $x_i = a_1 + id_1$ and $y_i = a_2 +
  id_2$ for $j = 0, 1, \dots, k-1$ are arithmetic progressions, and $\sigma$ is a permutation of $0,1,\ldots, k-1$.
\end{dfn}
Note that the appearance of this permutation is quite natural (and
necessary), since points with both coordinates in arithmetic
progression would all lie on a line. Garcia-Selfa and Tornero gave
examples of elliptic curves over $\mathbb{Q}$ that contain an
s.a.p.~of length 6. They also showed that there are only finitely many
such curves, and there are none with an s.a.p.~of length 7. Extending
their methods to s.a.p.'s of length 8 did not seem computationally
feasible, and they were not able to find an elliptic curve with an
s.a.p.~of length 8, or prove that none exists. The final open problem
they suggested is finding a universal bound for the length of
s.a.p.'s on elliptic curves over $\mathbb{Q}$.

In Section 2 of this paper we prove that 4319 is an upper bound for
the length of an s.a.p.~on an elliptic curve over $\mathbb{R}$, using
a combinatorial approach.  This solves the open problem above.  Given
a curve with a large s.a.p., we construct a graph on translates of the
curve, with many edges (segments of the curves) but not too many
vertices (translates of s.a.p.~points). Then we apply the well-known
crossing inequality to get a lower bound on the number of
intersections in the graph, and compare this with the upper bound that
we get from the fact that these low-degree curves cannot intersect too
often.

In Section 3 we generalize this method to arbitrary real algebraic
curves (not containing a line). We also generalize it from an
s.a.p.~to any $k^{2/3}$ points from a cartesian product of two length
$k$ arithmetic progressions. Then we give a second proof, using an old
result of Jarn\'{\i}k \cite{Jarn26}.  Finally, the result is extended to
complex algebraic curves (not containing a line).

\section{Elliptic curves over $\mathbb{R}$}

In this section we give a universal bound on the size of an s.a.p.~on
a real elliptic curve. We will use the following result a number of times.
\begin{lem}
\label{lem:intersect}
An elliptic curve (over $\mathbb{C}$) and a translate of that curve
can intersect in at most $4$ points (excluding points at infinity.)
\end{lem}
\begin{proof}
  Suppose the curve is given by $y^2+axy+by=x^3+cx^2+dx+e$.  A
  translate is given by
  $(y+v)^2+a(x+u)(y+v)+b(y+v)=(x+u)^3+c(x+u)^2+d(x+u)+e$ where at
  least one of $u,v$ does not equal $0$.  If $(x,y)$ is an
  intersection point of these curves then subtracting the one equation
  from the other we
  get \[2vy+v^2+avx+auy+auv+bv=3ux^2+3u^2x+u^3+2cux+cu^2+du.\] If
  $2v+au=0$ then all terms involving $y$ disappear.  In this case we
  have a quadratic in $x$ which can have at most $2$ real roots.
  Putting these values into the original equation we get at most $4$
  intersection points.  If $2v+au\ne 0$ then we can solve for $y$ to
  get \[y=\frac{3ux^2+(3u^2+2cu-av)x+(u^3+cu^2+du-auv-bv-v^2)}{2v+au}.\]
  Substituting this into the original equation we get $f(x)=0$ where
  $f$ is a quartic polynomial in $x$.  This polynomial has at most $4$
  roots.  Thus we cannot have more than $4$ intersection points of our
  elliptic curve and its translate.
\end{proof}

The main result is:
\begin{thm}
\label{thm:main}
  Consider an elliptic curve over a subfield of $\mathbb{R}$ given by
  $y^2+axy+by=x^3+cx^2+dx+e$.  Suppose we have an s.a.p.~on this curve
  of length $k$.  Then $k\le 4319$.
\end{thm}
The idea behind the proof is to consider translates of the curve in a
grid as a graph with edges between points in arithmetic progression
which occur consecutively on a translate.  We give a simple upper
bound on the crossing number of this graph and use the ``crossing
inequality'' to give a lower bound.  Putting these together we get the
stated upper bound for $k$.

\begin{dfn}
  Given a simple graph $G$, the crossing number, $\crn(G)$, is the
  minimum number of pairs of crossing edges in a planar drawing of
  $G$.
\end{dfn}

The crossing inequality was first proved independently by Ajtai, Chv\'atal, Newborn and Szemer\'edi \cite{Ajta82} and by Leighton \cite{Leig84}.  The version with the best bound to date, presented below, was given by Pach and T\'oth \cite{Pach97}.  \begin{thm}[Crossing Inequality] \label{thm:crossOrig}
  Suppose $G$ is a simple graph with $n$ vertices and $e$ edges.  If
  $e>7.5n$ then \[\crn(G)\ge\frac{e^3}{33.75n^2}.\]
\end{thm}

Pach and T\'oth also gave a crossing inequality for multigraphs which is the result we use herein.
\begin{thm}[Crossing inequality for multigraphs] \label{thm:cross}
  Suppose $G$ is a multigraph with $n$ vertices and $e$ edges (counting multiplicity.)  Suppose there are at most $m$ edges between any pair of vertices in $G$.  If
  $e>7.5mn$ then \[\crn(G)\ge\frac{e^3}{33.75mn^2}.\] \end{thm}

We have a few cases to consider depending on the number of connected
components of the curve and the number of points from the s.a.p.~on
each component.

\begin{figure}
\begin{center}
\subfigure[Two connected components]{
\includegraphics[width=57mm]{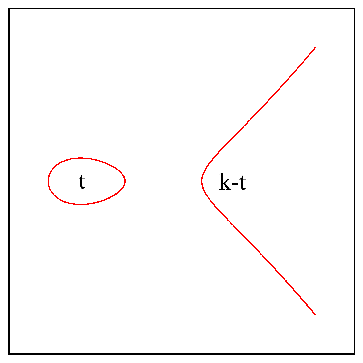}
\label{fig:twoComp} }
\subfigure[One connected component]{
\includegraphics[width=57mm]{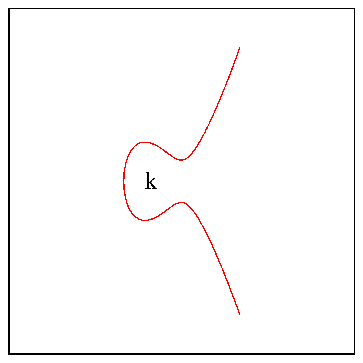}
\label{fig:oneComp} }
\end{center}
\caption{Elliptic curves over $\mathbb{R}$}
\label{fig:curves}
\end{figure}

Consider Figure \ref{fig:curves}.  First, suppose our curve has two
connected components with $t$ and $k-t$ points, respectively, of the
s.a.p.~on the components shown.  If there is more than one point on
the left component, so $t>1$, then we consider a graph containing
these $t$ points as vertices and the parts of the curve connecting
consecutive points as edges.  Then we clearly have a $t$-cycle.  If
there is only one point on the left component, so $t=1$, then consider
the graph with the vertex given by this point and no edges.  If $t=0$
then we only consider the connected component containing the point,
$[0:1:0]$, at infinity.  This case will be treated in the same way as
in the case where we only have one connected component.  So we need
only consider elliptic curves with two connected components.

\begin{figure}
\begin{center}
\includegraphics[width=57mm]{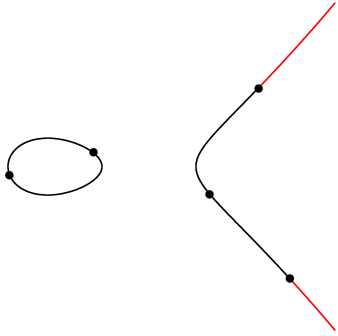}
\end{center}
\caption{Graph defined by an s.a.p.~on an elliptic curve}
\label{fig:graph}
\end{figure}

Consider the $k-t$ points on the component containing points at
infinity.  We extend the graph described above.  The idea is to
connect consecutive points, when considered as vertices, along the
curve with edges.  Then we will have $k-t$ vertices and $k-t-1$ edges.
Figure~\ref{fig:graph} gives an example of such a graph where $t=2$
and $k=5$.

We include the point $[0:1:0]$ in our graph to increase the number of
edges.  This is to improve the bound for $k$ in Theorem
\ref{thm:main}.  Connect the rightmost point on the top part of the
curve to infinity and do the same for the rightmost point on the
bottom part of the curve.  Our graph now contains $k-t+1$ vertices and
$k-t+1$ edges.  Considering the other component as well we end up with
$k+1$ vertices and $k+1$ edges if $t \ne 1$ and $k$ edges if $t=1$.
We will use variations of this graph in the proof of Theorem
\ref{thm:main}.

\begin{proof}[Proof of Theorem \ref{thm:main}]
  Suppose the s.a.p.~is given by $(x, y+\sigma(0)d_2), (x+d_1,
  y+\sigma(1)d_2), (x+2d_1, y+\sigma(2)d_2), \dots, (x+(k-1)d_1,
  y+\sigma(k-1)d_2)$.  Consider the $k^2$ translates of the elliptic
  curve given by all combinations of translating $x$ by $0, -d_1,
  -2d_1, \dots, -(k-1)d_1$ and $y$ by $0, -d_2, -2d_2, \dots,
  -(k-1)d_2$.  Considering the graph structure described above on
  these translates we get a graph with vertices given by all points in
  a $2k\times 2k$ grid and the point at infinity.  Thus we have
  $4k^2+1$ vertices.  We have to change the edges slightly to ensure
  that we have a well-defined graph.

  First note that we may have more than one edge connecting two
  vertices.  We show that the maximum multiplicity for such an edge is
  $4$.  If we have more than one edge connecting two vertices then
  these two vertices appear as consecutive points on a number of
  translates---see Figure~\ref{fig:edgeOne} for an example.  These
  points are given as $(x+ld_1, y+md_2)$ and $(x+l'd_1, y+m'd_2)$ for
  some $l, l', m, m'\in\mathbb{Z}$.  If these points appear on $r$
  translates then the difference vector $((l-l')d_1,(m-m')d_2)$
  connects $r$ pairs of points on the original elliptic curve---see
  Figure~\ref{fig:edgeTwo}.  But this is equivalent to having $r$
  points on the original curve intersecting $r$ points on a
  translate---see Figure~\ref{fig:edgeThree}.  Thus, by Lemma
  \ref{lem:intersect}, $r\le 4$.

\begin{figure}
\begin{center}
\subfigure[Three edges between a pair of vertices]{
\includegraphics[width=57mm]{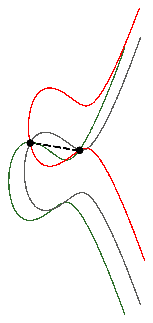}
\label{fig:edgeOne} }
\subfigure[Difference vector appearing on curve three times]{
\includegraphics[width=57mm]{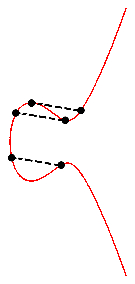}
\label{fig:edgeTwo} }
\subfigure[A pair of translates intersecting in three points]{
\includegraphics[width=57mm]{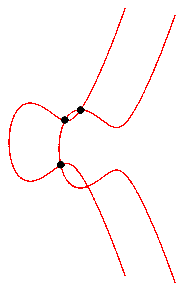}
\label{fig:edgeThree} }
\end{center}
\caption{The multiplicity of edges in the graph}
\label{fig:edgeCount}
\end{figure}

  Now since we are considering translates of a curve in a grid, a
  vertex may be a point on a number of curves.  Suppose $v_1$ and
  $v_2$ are consecutive points on a translate.  We may have a point
  $v_3$ on another translate which is actually between $v_1$ and $v_2$
  on the first translate.  In this case the edge from $v_1$ to $v_2$
  passes through the vertex $v_3$.  This is not allowed in a graph so
  we have to alter our graph slightly.  In this case we remove the
  edge in consideration from $v_1$ to $v_2$ and add an edge from $v_1$
  to $v_3$.  Performing this change where necessary we end up with a
  graph with the same number of vertices and edges but without the
  problem of an edge passing through a vertex to which it is not
  adjacent.  We call this graph $G$.

The only way we can have more than one edge going from a point in the
grid to the point at infinity is if that point is the rightmost point
on the top half of one translate and the rightmost point on the bottom
half of another translate.  Thus these edges have multiplicity at most
$2$.

Suppose $t$ is defined as in Figure \ref{fig:twoComp}.  If $t\ne 1$
then the number of edges, counting multiplicity, in $G$ is
$k^2(k+1)$, while if $t=1$ then the number of edges is $k^3$. We need 
only consider the case with less edges, so we assume we have $k^3$ edges. The number of vertices 
is $n=4k^2+1$.  Between any two
vertices there are at most $4$ edges.  Thus from the crossing inequality we
get \[\frac{(k^3)^3}{4(33.75)(4k^2+1)^2}\le \crn(G).\]

Any pair of translates intersect in at most $4$ points in the grid and
there are $\binom{k^2}{2}$ such pairs.  Thus the crossing number is
bounded by \[\crn(G)\le 4\binom{k^2}{2}.\]

When $t\ne 1$, putting these two inequalities together we
get \[\frac{(k^3)^3}{4(33.75)(4k^2+1)^2}\le
4\binom{k^2}{2}.\] Solving for $k$ in this inequality and noting that
$k$ is a positive integer we get $k\le 4319$.

Thus we have our uniform bound on the length, $k$, of an s.a.p.~on the
elliptic curve $y^2+axy+by=x^3+cx^2+dx+e$.
Note that in Theorem \ref{thm:cross} we require $e>4(7.5)n=30n$.  When
$t\ne 1$ this gives \[e>30n\Rightarrow
k^2(k+1)>30(4k^2+1)\Rightarrow k>15.\] When $t=1$ we have $k>16$.
\end{proof}

\section{General polynomials over $\mathbb{R}$}

In this section we generalize the result from the previous section to
arbitrary plane algebraic curves over $\mathbb{R}$. Above we did not 
fully use the structure of an s.a.p.; therefore we can also
generalize from an s.a.p.~to any $k^\alpha$ points from a $k\times k$ 
grid, where from the proofs below
we will see that we can take $\alpha = 2/3$.

We will give two proofs of this fact.  The first is a generalization of the proof used
above.  The second proof relies on a result of Vojt\v{e}ch Jarn\'{\i}k
\cite{Jarn26} about the possible number of lattice points on a convex
curve.  At the end of the section we show that the result for real
curves can be used to prove the result for complex curves.

By a $k\times k$ {\it grid} we mean the cartesian product of two
arithmetic progressions of length $k$; so an s.a.p.~consists of $k$
elements from a $k\times k$ grid with exactly one element on each row
and on each column.

The main result is:
\begin{thm} \label{thm:main2}
For every integer $d\geq2$ there is a constant $C$ depending only on $d$ such 
that if $f$ is a real plane algebraic curve of degree $d$ with no linear factor, 
then $f$ does not contain more than $C k^{\frac{2}{3}}$ points from a $k\times k$ grid.
\end{thm}

Both proofs use B\'ezout's Theorem.  For details see
\cite{Fult69}.
\begin{thm}[B\'ezout]
Suppose $F$ and $G$ are projective plane curves of degree $m$ and $n$
respectively defined over an algebraically closed field.  If $F$ and
$G$ do not have a common factor then they intersect in $mn$ points
counting multiplicity.
\end{thm}

Suppose $f(x,y)=0$ is an irreducible plane algebraic curve over
$\mathbb{C}$ of degree $d$.  If $f$ and a translate of $f$ have no
factors in common then homogenizing $f$ and the translate we can apply
B\'ezout's Theorem to get that these curves intersect in at most $d^2$
points.  Considering $f$ over a subfield, such as $\mathbb{R}$ or
$\mathbb{Q}$, this bound still holds.

\begin{rmks}
The number of multiple points on an irreducible plane curve $f$ of
degree $d$ is at most $(d-1)(d-2)/2$.\\
A result of Harnack gives that the number of connected components of
a real irreducible curve is at most $(d-1)(d-2)/2+1$. For details see
\cite{Gudk74}.
\end{rmks}

The first proof is almost identical to the proof for elliptic curves.
We consider the factor of $f$ that has the most points from the grid 
and construct a graph out of its $k^2$ translates on a $2k\times 2k$ grid.
A bound on the edge multiplicity is given by B\'ezout's Theorem. We now
have to deal with crossings given by self-intersections on our curve, as
well as with the possibility of many connected components. Fortunately,
by the remarks above, both of these are bounded by functions of the degree
of $f$.

\begin{proof}[First proof]
We will first assume that $f$ is irreducible. Suppose that $f$
contains $K=Ck^\alpha$ points from a $k \times k$ grid (we will
establish the appropriate value for $\alpha$ at the end).\\ We need to
assure that we do not have too many components with only one grid
point on them, since those points would not give any edges in the
graph.  A component with $m\ge 2$ grid points on it will give us $m$
edges, except that a component containing a point at infinity does not
form a closed loop, and only gives $m-1$ edges. But the line at
infinity and our curve can have at most $d$ intersections, by
Bez\'out's Theorem, thus we have at most $d$ such components.  By
choosing $C$ large enough, we can assure that we have at least $K/2$
(say) grid points that lie on a component not containing the point at
infinity.  This provides us with $K/2$ edges from each translate,
hence $k^2K/2$ edges in our graph.\\ The crossing number of the graph
is the number of intersections between translates plus the number of
self-intersections of translates.  By B\'ezout's Theorem, for any pair
of translates there are at most $d^2$ intersections. A
self-intersection is a multiple point, hence the number of these is
bounded by $(d-1)(d-2)/2$ as remarked above.  We get
\[cr(G)\le d^2\binom{k^2}{2}+k^2\frac{(d-1)(d-2)}{2}= \frac{1}{2}k^2(d^2k^2-3d+2)\le \frac{1}{2}d^2k^4.   \]
By the crossing inequality for multigraphs we get the lower
bound: 
\[\frac{(k^2(K/2))^3}{33.75d^2(4k^2)^2} \le \frac{e^3}{33.75d^2v^2}\le  cr(G).\]
Combining these we get $K^3 < C_1 d^4k^2$, with $C_1 = 2^6\cdot33.75$. Then going back to $K = Ck^{\alpha}$ gives
\[C^3k^{3\alpha-2} < C_1 d^4.\]
Now for $\alpha \ge 2/3$ we get a contradiction if we choose $C$ large
enough (depending only on $d$).  This proves the theorem for
irreducible curves, with constant $C = \sqrt[3]{C_1d^4}$.\\ For a
reducible curve $f$ with $Ck^{2/3}$ points from a $k\times k$ grid, we
have a factorization $f=f_1^{\alpha_1}f_2^{\alpha_2}\dots
f_r^{\alpha_r}$ where each $f_i$ is irreducible of degree $d_i$. Since
$f$ has no linear factor, we have $d_i \ge 2$ for all $i$, as well as
$r\le d/2$.  We take the factor $f_j$ which has the most points from
the grid on it, which is at least $\frac{C}{r}k^{2/3}\ge
\frac{2C}{d}k^{2/3}$.  Then by the result for irreducible curves, if
$C$ is large enough, we would get a contradiction. To be precise, we
need $\frac{2C}{d}\ge \sqrt[3]{C_1d^4}$, so $C \ge 13d^{7/3}$ would
do.
\end{proof}

For comparison, let's see the constant that we get this way when $f$
is an elliptic curve, so $d=3$.  Suppose the curve contains $k =
k^{1/3}\cdot k^{2/3}$ points of a $k\times k$ grid. Then we get a
contradiction when $k^{1/3}\geq 13d^{7/3} \geq 169$, or $k\geq
169^3\approx 5\cdot10^6$.\\ Using that the elliptic curve is
irreducible, we can instead take the inequality $K^3 < C_1d^4k^2$ from
the proof, with $K=k$.  Then the bound that we get is $k<C_1d^4 =
33.75\cdot 2^6\cdot 3^4 \approx2\cdot 10^5$.\\

The second proof uses a result of Jarn\'{\i}k \cite{Jarn26}.
\begin{thm} \label{thm:jarnik}
Suppose $f(x, y)$ is a strictly convex curve of length $N$.  Then the
number of integer points on $f$ is less than $cN^{2/3}$ for some
constant $c$.
\end{thm}

To apply this result we need to break up our curve into convex,
monotone pieces.  The number of such pieces is a function of the
degree of $f$.  Note that we actually break the curve up into convex
and concave pieces, but Theorem \ref{thm:jarnik} is valid for convex
or concave curves.  From now on when referring to a convex part of a
curve we will mean either a convex part or a concave part of the
curve.

We consider an irreducible curve $f$ of degree $d$.  To break $f$ up
into convex parts we need to cut the curve at all inflection points
and all singularities, i.e. points where either of the first derivatives vanish.

The following result bounds the number of inflection points.  For a
proof see \cite{Kirw92} or the exercises in \cite{Fult69}.
\begin{lem}
Suppose $f$ is an irreducible curve of degree $d$.  Then $f$ has at
most $3d(d-2)$ inflection points.
\end{lem}

Now we consider the points where $f_x=0$.  By assumption $f$ is
irreducible so $f$ and $f_x$ are coprime.  Also, the degree of $f_x <
d$.  Thus we can apply B\'ezout's Theorem to get that there are at most
$d(d-1)$ points where $f=0$ and $f_x=0$.  Similarly, there are at most
$d(d-1)$ points where $f=0$ and $f_y=0$.

So we need at most $3d(d-2) + 2d(d-1) = d(5d-8)$ cuts to break $f$
into convex, monotone parts.

\begin{proof}[Second proof (of Theorem \ref{thm:main2})]
Suppose again that $f$ is irreducible and contains $K = Ck^\alpha$
points from a $k\times k$ grid.  Firstly we scale and translate $f$ so
that the gap in the $k\times k$ grid is $1$ in both the $x$- and
$y$-directions and the points of the grid are integral.  Convexity is
preserved under this transformation.  Now we can separate $f$ into
convex, monotone parts using at most $d(5d-8)$ cuts.  One of these
parts has at least the average number of points from the grid on $f$.
So we have at least $K/d(5d-8)$ points from the grid on this part of
the curve.  Since the grid has gap $1$ and length $k$ we can bound the
length of this part of the curve by $2k$.  Thus we get, by Theorem
\ref{thm:jarnik}, that
\[\frac{K}{d(5d-8)} < c(2k)^{2/3}.\]
This gives
\[C k^{\alpha-2/3} < c_1 d^5,\]
so again we get a contradiction for $\alpha \ge 2/3$ and $C$ large enough.

Using the method at the end of the first proof of
Theorem~\ref{thm:main2} we get the result for reducible curves.
\end{proof}

The dependence on $d$ of the constant $C$ in the theorem cannot be removed, 
as the following example shows.

\begin{rmk}
Given $d(d+3)/2$ points in the plane there is a curve of degree $d$
passing through those points.
\end{rmk}
\begin{proof}
An arbitrary curve $f(x,y)$ of degree $d$ contains $(d+1)(d+2)/2$
terms.  To see this note that there are $d+1$ monomials in $x$ and $y$
of degree $d$, $d$ monomials of degree $d-1$ and so on.  But we are
considering $f(x,y)=0$ and so one of the terms is dependent on the
others.  Thus we have $(d+1)(d+2)/2-1=d(d+3)/2$ terms in $f$.

Now, given $d(d+3)/2$ points in the plane, we can plug each of these
into $f$.  This gives us $d(d+3)/2$ linear equations in $d(d+3)/2$
unknowns.  Thus a solution exists and so we can find a curve of degree
$d$ going through the $d(d+3)/2$ points.
\end{proof}

Consider any s.a.p.~of length $k$.  By the above remark there exists a
curve of degree $d$ containing the s.a.p.~where $k=d(d+3)/2$ and so
$d=(-3+\sqrt{9+8k})/2$.\\

Theorem \ref{thm:main2} can be extended to any complex algebraic plane
curve.  We use the result for the reals to prove the result for the
complex case. By a $k\times k$ grid in the complex plane we shall mean
a cartesian product of two arithmetic progressions in the complex plane. 
By an arithmetic progression in the complex plane we mean points 
$\alpha+i\beta$ with $\alpha,\beta\in \mathbb{C}$ and $i = 0,1,\ldots,k-1$.

\begin{thm}
For every integer $d\geq2$ there is a constant $C$ depending only on $d$ such 
that if $f$ is a \emph{complex} plane algebraic curve of degree $d$ with no linear factor, 
then $f$ does not contain more than $C k^{\frac{2}{3}}$ points from a $k\times k$ grid.
\end{thm}
\begin{proof}
Suppose $f(w,z)$ is our complex curve with many points on a $k\times
k$ grid, given by $\alpha + j\beta$ in one direction and $\gamma +
j\delta$ in the other direction where $\alpha, \beta, \gamma, \delta
\in \mathbb{C}$ and $j=0, 1, \dots, k-1$.  Now consider the curve
$g(x, y)=f(\alpha + x\beta, \gamma + y\delta)$.  This is a curve with
complex coefficients in two real variables and many points on the
$k\times k$ grid consisting of the points $(i,j)$ with $i,j =
0,1,\ldots,k-1$. The real and imaginary parts of $g$ are real
algebraic curves each having many points on the same grid.  Thus by
Theorem \ref{thm:main2} there is a large constant $C$ for which we get
a contradiction.
\end{proof}

Jarn\'{\i}k's result gives the existence of a bound in
Theorem~\ref{thm:main2}.  This bound is by no means optimal.  In fact,
Bombieri and Pila proved in \cite{Bomb89} that we can get the bound
$c(d,\varepsilon)k^{1/d+\varepsilon}$ for any $\varepsilon>0$ if the
curve is irreducible.  This clearly gives a better bound for large
degree.

\section*{Acknowledgment}

The authors would like to thank Trevor Wooley for his helpful remarks.

\bibliographystyle{plain} \bibliography{references}

\begin{thebibliography}{10}

\bibitem{Ajta82}
M.~Ajtai, V.~Chv{\'a}tal, M.~M. Newborn, and E.~Szemer{\'e}di.
\newblock Crossing-free subgraphs.
\newblock In {\em Theory and practice of combinatorics}, volume~60 of {\em
  North-Holland Math. Stud.}, pages 9--12. North-Holland, Amsterdam, 1982.

\bibitem{Bomb89}
E.~Bombieri and J.~Pila.
\newblock The number of integral points on arcs and ovals.
\newblock {\em Duke Mathematical Journal}, 59(2):337--357, 1989.

\bibitem{Brem99b}
A.~Bremner.
\newblock On arithmetic progressions on elliptic curves.
\newblock {\em Experimental Mathematics}, 8(4):409--413, 1999.

\bibitem{Brem99a}
A.~Bremner.
\newblock On squares of squares.
\newblock {\em Acta Arith.}, 88(3):289--297, 1999.

\bibitem{BrST00}
A.~Bremner, J.~H. Silverman, and N.~Tzanakis.
\newblock Integral points in arithmetic progression on $y^2 = x(x^2 − n^2)$.
\newblock {\em J. Number Theory}, 80:187--208, 2000.

\bibitem{Fult69}
William Fulton.
\newblock {\em Algebraic Curves: {A}n {I}ntroduction to {A}lgebraic
  {G}eometry}.
\newblock 1969.

\bibitem{Garc06}
I.~Garcia-Selfa and J.M. Tornero.
\newblock On simultaneous arithmetic progressions on elliptic curves.
\newblock {\em Experimental Mathematics}, 15(4):471--478, 2006.

\bibitem{Gudk74}
D.A. Gudkov.
\newblock The topology of real projective algebraic varieties.
\newblock {\em Russian Mathematical Surveys}, 29(4):1--79, 1974.

\bibitem{Jarn26}
V.~Jarn\'{\i}k.
\newblock \"{U}ber die {G}itterpunkte auf konvexen {K}urven.
\newblock {\em Mathematische Zeitschrift}, 24(1):500--518, 1926.

\bibitem{Kirw92}
F.C. Kirwan.
\newblock {\em Complex Algebraic Curves}.
\newblock Cambridge University Press, 1992.

\bibitem{Leig84}
Frank~Thomson Leighton.
\newblock New lower bound techniques for {VLSI}.
\newblock {\em Math. Systems Theory}, 17(1):47--70, 1984.

\bibitem{Pach97}
J\'{a}nos Pach and G\'{e}za T\'{o}th.
\newblock Graphs drawn with few crossings per edge.
\newblock In {\em GD '96: Proceedings of the Symposium on Graph Drawing}, pages
  345--354, London, UK, 1997. Springer-Verlag.

\end{thebibliography}
\end{document}